\documentclass[12pt]{amsart}
\usepackage{amsthm,amssymb}
\input epsf.tex

\def\crossout#1{\backslash\!\!\!/\!\!\!#1}
\def\disk{\mathbb D}
\def\euclid{\mathbb E}
\def\bbK{\mathbb K}
\def\Hyp{\mathbb H}
\def\R{\mathbb R}
\def\sphere{\mathbb S}
\def\dsph{d_{\sphere^n}}

\def\norm#1{\left\|#1\right\|}

\def\sH{\mathcal H}
\def\sE{\mathcal E}
\def\sS{\mathcal S}
\def\sT{\mathcal T}
\def\sZ{\mathcal Z}

\newtheorem{prop}{Proposition}
\newtheorem{thm}[prop]{Theorem}
\newtheorem{lemma}[prop]{Lemma}
\newtheorem{coro}[prop]{Corollary}

\theoremstyle{definition}
\newtheorem{rmk}{Remark}

\title[]{Comparing Corresponding Dihedral Angles on Classical Geometric Simplices}
\author[Thomas Au]{Thomas Kwok-keung Au}
\address{Department of Mathematics, The Chinese University of Hong
Kong, Shatin, Hong Kong.} \email{thomasau@cuhk.edu.hk}
\author[Feng Luo]{Feng Luo}
\address{Center of Mathematical Sciences, Zhejiang University,
Hangzhou, China. \vspace*{-3ex}} \address{Department of Mathematics,
Rutgers University, Piscataway, NJ08854, USA}
\email{fluo@math.rutgers.edu}
\author[Richard Stong]{Richard Stong}
\address{Department of Mathematics, Rice University, Houston,
TX77005, USA} \email{stong@math.rice.edu}

\begin{document}
\maketitle
\begin{abstract}
In this article, we prove a theorem comparing the dihedral angles of
simplexes in the  hyperbolic, spherical and Euclidean geometries.
\end{abstract}

\setlength\baselineskip{24pt} \setlength\parindent{0em}
\setlength\parskip{1.5ex} \thispagestyle{empty}

\section{Introduction}
It is well known that given any spherical (or hyperbolic) triangle,
one can decrease (or increase) its inner angles to obtain an
Euclidean triangle. The goal of this paper is to establish this
general fact for all dimensions. It was motivated by the study of
the volume of convex polytopes in classical geometry in terms of
dihedral angles.  Interesting related topics can be found in
\cite{Luo2,Luo3}.

By a space of classical geometry, we mean the
$n$-sphere~$\sphere^n$, the Euclidean $n$-space~$\euclid^n$, or the
hyperbolic $n$-space~$\Hyp^n$. For simplicity, they will be
collectively denoted by $\bbK^n$. A classical geometric $n$-simplex
$\sZ$ in $\bbK^n$ or simply a $\bbK^n$-simplex is the geodesic
convex hull of $(n+1)$ points $z_1,$ $z_2, \ldots$, $z_{n+1}$ in
$\bbK^n$ so that these points are not lying in any
$(n-1)$-dimensional totally geodesic submanifold. These $n+1$~points
are called the vertices of the simplex~$\sZ$.   As a convention, we
always consider simplexes with vertex labelled. That is, the simplex
$\sZ$ is represented by the $(n+1)$-tuple $( z_1, z_2, \ldots,
z_{n+1}) \in \left(\bbK^n\right)^{n+1}$, where $z_i$ denotes the
$i$-th vertex. In addition, two simplexes are equivalent if there is
a $\bbK^n$-isometry taking such an $(n+1)$-tuple to another.

We will compare the dihedral angles of simplexes in classical
geometries. Let $\sZ = ( z_1, z_2, \ldots, z_{n+1})$ be a
$\bbK^n$-simplex. We denote the codimension-1 face opposite to the
$i$-th vertex $z_i$ by $F_i(\sZ)=( z_1, ..., ,\crossout{z}_i, ...,
z_{n+1}) \in \left(\bbK^n\right)^n$. Then the dihedral angle
$\zeta_{ij}$, for $i\ne j$, is the angle between the faces
$F_i(\sZ)$ and $F_j(\sZ)$. Let $\sS$ and $\sT$ be two
$\bbK^n$-simplexes (not necessarily the same $\bbK^n$) of dihedral
angles $\sigma_{ij}$ and $\tau_{ij}$ respectively.  It is said that
$\sS \preceq \sT$ if and only if $\sigma_{ij} \leq \tau_{ij}$ for
every $i, j$. If in addition, there is a pair of $i\ne j$ such that
$\sigma_{ij} < \tau_{ij}$, then $\sS \prec \sT$.

In this article, we are going to demonstrate a theorem which may be
roughly abbreviated by $\Hyp^n \prec \euclid^n \prec \sphere^n\,.$
\renewcommand\theprop{(Comparison of Simplexes)}
\begin{thm}
There is a natural partial order on $n$-simplexes in these spaces of
classical geometry according to dihedral angles. More precisely,
\begin{description}
\item[M1] For every $\sphere^n$-simplex $\sS$, there is an $\euclid^n$-simplex
$\sE$ such that $\sE\prec \sS$.

\item[M2] For every $\Hyp^n$-simplex $\sH$, there is an $\euclid^n$-simplex
$\sE$ such that $\sH\prec \sE$.

\item[M3] For every $\euclid^n$-simplex $\sE$, there is an $\sphere^n$-simplex
$\sS$ and an $\Hyp^n$-simplex $\sH$ such that $\sH \prec \sE \prec
\sS$.

\item[M4] if $\sE_1$, $\sE_2$ are $\euclid^n$-simplexes such that $\sE_1 \preceq \sE_2$,
then $\sE_1$ and $\sE_2$ have exactly the same corresponding
dihedral angles.

\end{description}
\end{thm}
\renewcommand\theprop{\arabic{prop}}
\addtocounter{prop}{-1}
\begin{rmk}
The statement~{\bf M4} above was also proved by Richard Stong.
Moreover, the theorem is trivial for $n=2$.
\end{rmk}
The statements {\bf M3} and {\bf M4} are proved in \S 2 by
considering suitable variations of Gram matrices.  In a certain
sense, Euclidean Gram matrices lie in the common boundary of
spherical and hyperbolic ones.  In \S 3, we will prove {\bf M2}
using geometric comparison.  From a given $\Hyp^n$-simplex, the
desired $\euclid^n$-simplex has the same inscribed sphere and the
dihedral angle comparison follows naturally from Gauss-Bonnet
Theorem. To prove {\bf M1}, on the one hand, the method of Gram
matrices fails because we do not have control of the signs of
cofactors. On the other hand, the geometric construction of a
compact $\euclid^n$-simplex from a given $\sphere^n$-simplex is more
subtle. The idea is to ``extend'' or ``enlarge'' the dual of the
given $\sphere^n$-simplex. Then take the Euclidean dual of the
``extended'' simplex and perturb a little bit if necessary. The
details will be discussed in \S 4.

\section{Gram Matrices}
The Gram matrix $G=G(\sZ)$ of a $\bbK^n$-simplex~$\sZ$ is an
$(n+1)\times(n+1)$ matrix with entries $-\cos\zeta_{ij}$, where
$\zeta_{ij}$ is the dihedral angles of $\sZ$ with the convention
$\zeta_{ii}=\pi$. It is clearly symmetric and has diagonal entries
equal to~1.  Since the function $-\cos(\,\cdot)$ is monotonic
increasing on $(0,\pi)$, it is also natural to say that two Gram
matrices $(a_{ij}) \preceq (b_{ij})$ if their corresponding entries
$a_{ij} \leq b_{ij}$ for all $i, j$.

First, let us recall a result of \cite{Luo} and \cite{Milnor} which
clarifies the relation between Gram matrices and classical geometric
simplexes.
\begin{thm}\label{thm-GramSimplex}
Let $A$ be an $(n+1)\times(n+1)$ real symmetric matrix with diagonal
entries equal~$1$ and let $c_{ij}$ be the $(i,j)^{\text{th}}$
cofactor of~$A$.
\begin{enumerate}
\item $A$ is the Gram matrix of an $\sphere^n$-simplex if and only if $A$ is
positive definite.

\item $A$ is the Gram matrix of an $\euclid^n$-simplex if and only if
$\det(A)=0$,  all principal $n\times n$ submatrices of~$A$ are
positive definite, and all $c_{ij} > 0$.

\item $A$ is the Gram matrix of an $\Hyp^n$-simplex if and only if
$\det(A)<0$, all principal $n\times n$ submatrices of~$A$ are
positive definite, and all $c_{ij}>0$.
\end{enumerate}
\end{thm}
For simplicity, we may refer to the above cases of Gram matrices as
spherical, Euclidean, or hyperbolic Gram matrices.

Using continuous variation of Gram matrices, we are able to show
that an Euclidean simplex sits between a hyperbolic and a spherical
ones.
\renewcommand\theprop{M3}
\begin{thm}
For any Euclidean $n$-simplex $\sE$, there is a hyperbolic
$n$-simplex $\sH$ and a spherical $n$-simplex $\sS$ such that $\sH
\prec \sE \prec \sS$.
\end{thm}
\renewcommand\theprop{\arabic{prop}}\addtocounter{prop}{-1}
\begin{proof}
Let $\sE$ be an $\euclid^n$-simplex and $G=\left(g_{ij}\right)$ be
its corresponding Gram matrix with cofactors $c_{ij}$. In other
words, by Theorem~\ref{thm-GramSimplex}, $\det(G)=0$ and $c_{ij}>0$
for all $i, j$, together with all principle $n\times n$ submatrices
of $G$ being positive definite.  Let $P=\left(p_{ij}\right)$ be the
$(n+1)\times(n+1)$ matrix in which every diagonal entry is~$1$ and
$p_{ij}\equiv-1$ for all $i\ne j$. Let $A(t) =
\left(\,a_{ij}(t)\,\right)$ be the path in the space of
$(n+1)\times(n+1)$ symmetric matrices defined by,
$$
A(t) = (1-t)G + t\,P \qquad t\in [0,1].
$$
It is clear that the eigenvalues of the principal $n\times n$
matrices of $A(t)$ and the cofactors $c_{ij}(t)$ of $A(t)$ depend
continuously on the entries of~$A(t)$ and hence in~$t$.  Thus, for
sufficiently small~$t>0$, the principal $n\times n$ matrices remain
positive definite and $c_{ij}(t)>0$.  Moreover,
$$
\frac{d}{dt}\left[ \det A(t) \right] = \sum_{i,j=1}^{n+1}
c_{ij}(t)a_{ij}'(t) = \sum_{i\ne j} c_{ij}(t)(-1-g_{ij}).
$$
Since $c_{ij}(0)>0$, for sufficiently small $t>0$, we have $\det
A(t) <0$.  Thus, by Theorem~\ref{thm-GramSimplex}, $A(t)$
corresponds to the Gram matrix of a $\Hyp^n$-simplex~$\sH$. Clearly,
$a_{ij}(t) < g_{ij}$ for $i\ne j$.

To obtain an $\sphere^n$-simplex~$\sS$, one simply takes another
matrix~$P$ which has all entries $p_{ij}\equiv 1$ for all $i, j$.
This clearly produces $a_{ij}(t) > g_{ij}$ for $i\ne j$. The
argument is exactly the same as above with the only difference that
$\det(A(t))>0$. Again, $A(t)$ corresponds to the Gram matrix of a
$\sphere^n$-simplex~$\sS$. We then have $\sH \prec \sE \prec \sS$.
\end{proof}
\begin{rmk}
From the proof, we actually have $\sS$ and $\sH$ which have dihedral
angles arbitrarily close to those of $\sE$.
\end{rmk}

The Gram matrices also provides another proof for the ``rigidity''
of Euclidean simplexes given by Stong.
\renewcommand\theprop{M4}
\begin{thm}
If $\sE_1$ and $\sE_2$ are two Euclidean $n$-simplexes such that
$\sE_1 \preceq \sE_2$, then they are similar.
\end{thm}
\renewcommand\theprop{\arabic{prop}}\addtocounter{prop}{-1}
\begin{proof}
Let $\sE_1$ and $\sE_2$ be two Euclidean $n$-simplexes such that
$\sE_1 \preceq \sE_2$. Furthermore, let $G_1$ and $G_2$ be their
corresponding Gram matrices and $A(t) = (1-t)G_1 + tG_2$,
$t\in[0,1]$ be a path in symmetric matrices joining the two Gram
matrices.  We also denote the cofactors of $A(t)$ by $c_{ij}(t)$.

By Theorem~\ref{thm-GramSimplex}, all principal $n\times n$
submatrices of $G_1$ and $G_2$ are positive definite and $\det(G_1)
= 0 = \det(G_2)$.  Thus, both $G_1$ and $G_2$ are semi-positive
definite. As a consequence, $A(t)$ is semi-positive definite for
all~$t$.  In particular, $\det(A(t)) \geq 0$ for all $t\in[0,1]$.
Let $f(t) = \det(A(t))$.  It is obvious that
$$
f'(t) = \frac{d}{dt}\left( \det(A(t)) \right) = \sum_{i\ne j}
\left(\cos\alpha_{ij} - \cos\beta_{ij}\right) c_{ij}(t),
$$
where $\alpha_{ij}$ and $\beta_{ij}$ are the dihedral angles of the
simplexes $\sE_1$ and $\sE_2$ respectively.  Note that for all
$i,j$, $\alpha_{ij}\leq\beta_{ij}$, thus
$\cos\alpha_{ij}\geq\cos\beta_{ij}$.  Suppose there is a pair of
corresponding dihedral angles $\alpha_{pq} < \beta_{pq}$. Since
$G_1$ and $G_2$ are Euclidean Gram matrices, for all $i,j$, we have
$c_{ij}(0)> 0$ and $c_{ij}(1)
> 0$.  As a consequence,
$$
f'(1) \geq \left(\cos\alpha_{pq} - \cos\beta_{pq}\right) c_{pq}(1)
> 0.
$$
Together with the fact that $f(1) = 0$, there is a small
$\varepsilon > 0$ such that $f(t) < 0$ for $t\in(1-\varepsilon,1)$.
This contradicts that $f(t)\geq 0$. Hence, for all $i,j$, one must
have $\alpha_{ij} = \beta_{ij}$.
\end{proof}

\section{Gauss-Bonnet}
To show that a hyperbolic simplex is dominated by an Euclidean one,
it only requires a simple geometric construction and an angle
comparison based on the Gauss-Bonnet Theorem.
\renewcommand\theprop{M2}
\begin{thm}
For every hyperbolic $n$-simplex $\sH$, there is an Euclidean
$n$-simplex $\sE$ such that $\sH \prec \sE$.
\end{thm}
\renewcommand\theprop{\arabic{prop}}\addtocounter{prop}{-1}
Let $\sH$ be a $\Hyp^n$-simplex in the Poincar\'e disc model
$\disk^n$ of the hyperbolic space.  Let $S\subset\disk^n$ be an
inscribed hyperbolic $(n-1)$-sphere of~$\sH$. Without loss of
generality, by a hyperbolic isometry, one may assume that the
in-center of~$\sH$ is the origin and so $S$ is an Euclidean sphere
with center at the origin.

Let $u_1, \ldots, u_{n+1} \in\disk^n$ be the points of tangency
of~$S$ with~$\sH$.  They are also considered as Euclidean vectors
from the origin.  Let us first give an algebraic description of the
geometry of the vectors.
\begin{lemma}\label{lem-neversided}
\begin{enumerate}
\item Any $n$ vectors among $\left\{ u_1, \ldots, u_{n+1} \right\}$ are
linearly independent.

\item The system of linear equations $\displaystyle \sum_{i=1}^{n+1} x_i u_i =
0$ has only a 1-dimensional solution space of the form $\left( x_1,
\ldots, x_{n+1} \right)$ where $x_i x_j >0$ for all $i,j$.  That is,
the $x_i$'s are all of the same sign.
\end{enumerate}
\end{lemma}
\begin{proof}
Let $W_1, \ldots, W_{n+1} \subset \disk^n$ be codimension-1
hyperbolic hypersurfaces tangent to~$S$ at $u_1, \ldots, u_{n+1}$
respectively. That is, they determine the $(n-1)$-faces of~$\sH$.
Since $\sH$ is nondegenerate, the first statement is evident;
otherwise, there will be $n$~such faces intersecting in a
1-dimensional geodesic but not a vertex.

Suppose the second statement is not true.  By simple Linear Algebra,
there is a vector~$v$ satisfying $\langle v,u_i \rangle \geq 0$ for
all~$i=1,\ldots,n+1$.  Take a geodesic~$L$ from the center of~$S$
along the direction~$v$.  This geodesic makes an angle~$\geq\pi/2$
with~$u_i$ and so does not intersect any $W_i$.  Otherwise, there
will be a hyperbolic triangle with angle sum~$>\pi$.  Hence, the
hypersurfaces $W_i$ do not bound a compact simplex.
\end{proof}
\begin{rmk}
Note that in the disk model, the geodesic~$L$ from the center is
also an Euclidean ray.  Thus, the same argument proves an analogue
of the second statement in $\euclid^n$.
\end{rmk}
\begin{proof}[Proof of\/ {\bf M2}] Let $P_1,
\ldots, P_{n+1}$ be the Euclidean codimension-1 hyperplanes
in~$\R^n$ tangent to~$S$ at $u_1, \ldots, u_{n+1}$ respectively.
Then by Lemma~\ref{lem-neversided}, an $\euclid^n$-simplex~$\sE$ is
bounded by $P_1, \ldots, P_{n+1}$ with the origin as its in-center.

Since each $P_i$ has normal vector~$u_i$, the dihedral angles
$\xi_{ij}$ of $\sE$ are given by
$$
\xi_{ij} = \pi - \angle(u_i,u_j) = \pi - \arccos \dfrac{\langle u_i,
u_j\rangle} {\|u_i\|\,\|u_j\|}.
$$
As in the above lemma, we continue to use $W_1,\ldots,W_{n+1}$ to
denote codimension-1 hyperbolic hypersurfaces with normals $u_1,
\ldots, u_{n+1}$.  Then these $W_i$'s bound the hyperbolic simplex
$\sH$. Let $\eta_{ij}$ be the hyperbolic dihedral angle between
$W_i$ and $W_j$.  Consider the Euclidean 2-plane $P$ through the
origin spanned by the vectors $u_i$ and $u_j$. Then by the
construction $P$ is perpendicular to both $P_i$ and $P_j$.  Let
$D^2$ be the intersection $P \cap \disk^n$.  Then $D^2$  is a
totally geodesic hyperbolic 2-plane perpendicular to $W_i$ and $W_j$
(see the figure below).
\begin{center}
\mbox{\epsfysize=33mm \epsfbox{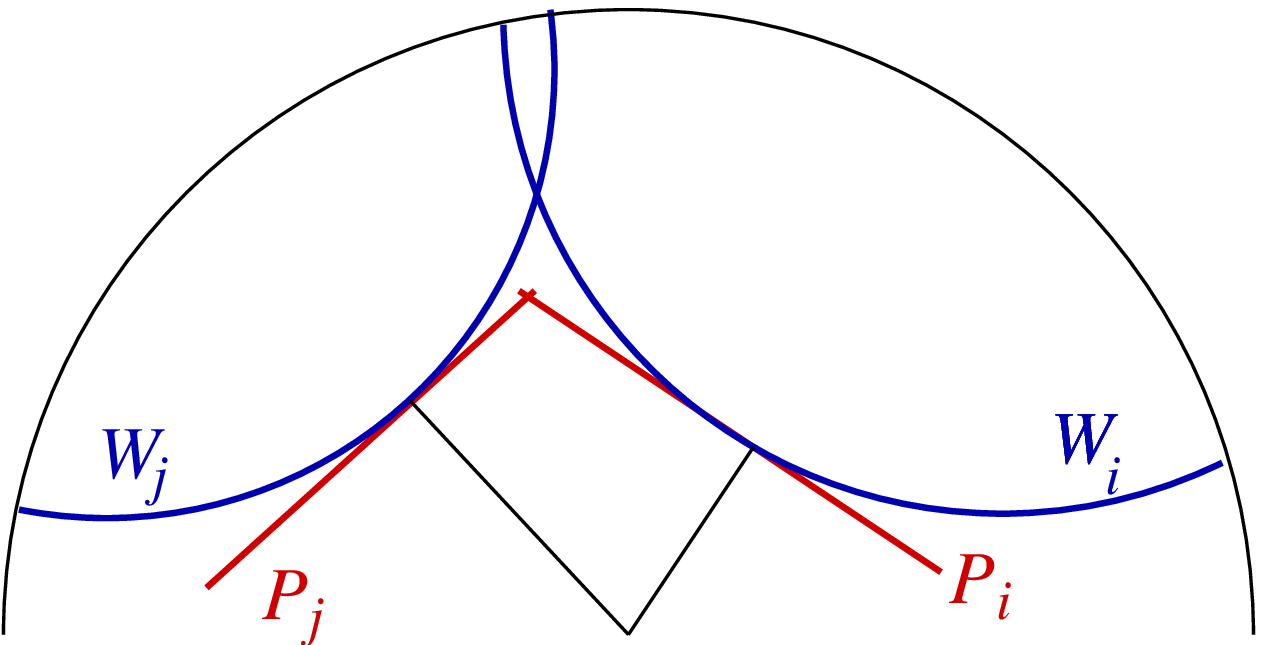}}
\end{center}
The intersections of $D^2$ with $W_i$ and $W_j$ respectively produce
two geodesics $\gamma_i$ and $\gamma_j$ in $D^2$. These two
geodesics together with the geodesics $u_i$ and $u_j$ from the
origin form a hyperbolic quadrilateral in $D^2$ with inner angles
$\pi-\xi_{ij}, \pi/2, \pi/2, \eta_{ij}$.  By Gauss-Bonnet Theorem,
it follows that their sum is less than~$\pi$.  Thus
$\eta_{ij}<\xi_{ij}$ and hence $\sH \prec \sE$.
\end{proof}

\section{The Sphere}
In this last section, we will deal with spherical simplexes. The
following convention will be adopted.  Let $\sphere^n$ be the unit
sphere in~$\euclid^{n+1}$; $\euclid^n = \euclid^n\times\{0\} \subset
\euclid^{n+1}$ and $\sphere^{n-1} = \sphere^n \cap \euclid^n$.
\renewcommand\theprop{M1}
\begin{thm}
For any  spherical $n$-simplex~$\sS$ with dihedral
angles~$\sigma_{ij}$, there is an Euclidean $n$-simplex~$\sE$ with
dihedral angles~$\xi_{ij}$ such that $\xi_{ij} < \sigma_{ij}$ for
all $i,j$.
\end{thm}
\renewcommand\theprop{\arabic{prop}}\addtocounter{prop}{-1}

The strategy of the proof goes as follows.  Let $\sS \subset
\sphere^n$ be a spherical simplex with dihedral angles
$\sigma_{ij}$, $i,j = 1,\ldots,n+1$.  Consider its dual
$\sphere^n$-simplex $\sS^*= \left( v_1,\ldots,v_{n+1}\right) \in
\left(\sphere^n\right)^{n+1}$.  By duality, the spherical distance
between the vertices is given by $\dsph(v_i,v_j) = \pi-\sigma_{ij}$.
We will move the vertices $v_i$'s appropriately to increase the
distances $\dsph(v_i, v_j)$ until it becomes the spherical dual of
an Euclidean $n$-simplex. The Euclidean $n$-simplex will have
dihehral angles smaller than $\sigma_{ij}$.

In  the rest of the section, for a $k$-ball~$B\subset\sphere^n$, by
\it an hemi-sphere \rm in $\partial B$, we refer to a closed
$(k-1)$-ball in $\partial B$ of the same radius as~$B$.

First, let us recall briefly the dual of an Euclidean $n$-simplex
$\sE$ in $\euclid^n$.  The following is a well-known fact.  See, for
instance, \cite{Luo} for a proof.
\begin{lemma}\label{lemma-EuclidCpt}
Given $n+1$ points $w_1, \ldots, w_{n+1} \in \sphere^{n-1} \subset
\euclid^n$, the convex polytope $\sE = \left\{ x\in\euclid^n :
\langle x-w_i, w_i\rangle \leq 0~\text{for all $i$}\right\}$ bounded
by the tangent planes to $\sphere^{n-1}$ at $w_i$'s in the side
containing the origin is an Euclidean $n$-simplex~$\sE$ if and only
if $\left\{ w_1, \ldots, w_{n+1} \right\}$ does not lie in any
hemi-sphere in $\sphere^{n-1}$.
\end{lemma}
We call $\left( w_1, \ldots, w_{n+1} \right) \in
\left(\sphere^{n-1}\right)^{n+1}$ the {\sl spherical dual\/}
of~$\sE$. Note that the $(i,j)^{\text{th}}$ dihedral angle of $\sE$
is $\pi-\dsph(w_i,w_j)$.

Second, we need a process of extending the sides of a geodesic
triangle on~$\sphere^n$.
\begin{lemma}\label{lemma-extend}
Let $T_0$ be a spherical triangle of angles $a$, $b$, $c$ and
corresponding opposite side lengths $x(0)$, $y(0)$, $z(0)$. Let
$T_t$ be a 1-parameter family of spherical triangles obtained by
extending the geodesics of lengths $x(0)$ and $y(0)$ to $x(t)$ and
$y(t)$ respectively in the same growth rate $x'(t)=y'(t)=g(t)>0$
while keeping the angle $c$ is fixed. If $x(t)+y(t)<\pi$, then the
length $z(t)$ of the third side satisfies $z(t) > z(0)$.
\end{lemma}
\begin{center}
\mbox{\epsfysize=38mm \epsfbox{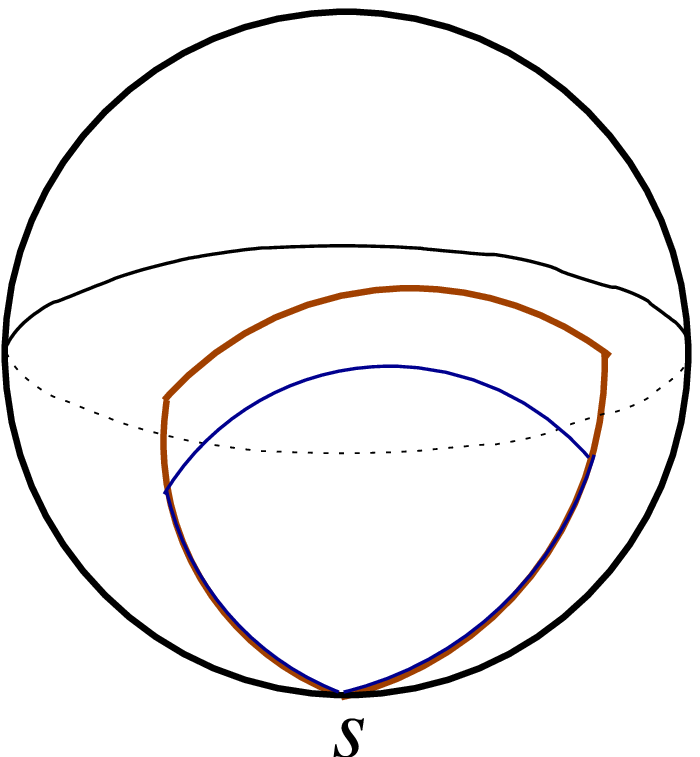}}
\end{center}
\begin{proof}
According to the spherical Cosine Law, for each $T_t$, we have
$$
\cos z(t) = \cos x(t)\cos y(t) + \sin x(t)\sin y(t)\,\cos(c).
$$
Differentiating with respect to~$t$ and grouping terms, we have
\begin{align*}
z'(t)\sin z(t) &=  x'(t)\sin x(t)\cos y(t) + y'(t)\cos x(t)\sin y(t) \\
&\qquad\quad {} - x'(t)\cos x(t)\sin y(t)\cos(c) - y'(t)\sin x(t)\cos y(t)\cos(c) \\
&= g(t)\,(1-\cos(c)) \sin( x(t)+y(t) ) > 0.
\end{align*}
Thus, $z(t)$ keeps increasing as long as the condition
$x(t)+y(t)<\pi$ holds.
\end{proof}
To begin the proof, consider the dual simplex $\sS^* =
\left(v_1,\ldots,v_{n+1}\right) \in \left(\sphere^n\right)^{n+1}$ of
the given one~$\sS$.  Let $B_s\subset\sphere^n$ be the spherical
$n$-ball of the smallest radius containing~$\sS^*$. Without loss of
generality, assume its center is located at~$s=(0,\ldots,0,-1)\in
 \bf \euclid^{n+1}$.  Evidently, its radius~$<\pi/2$ and there are at
least two vertices $v_i$'s lying on the boundary of $B_s$.  By
permutating the vertex labels, we may assume that $v_1, \ldots, v_m
\in
\partial B_s$ for $2\leq m \leq n+1$, while $v_{m+i}\in int(B_s)$ for $i\geq 1$.

For each vertex~$v_i\in\sS^*$, let $\gamma_i(t)$ be the unique
geodesic ray from~$s$ to~$-s$ through~$v_i$ such that $t \in
[-d_{\sphere^n}(s, v_i), \infty)$ is the arc length parameter with
$\gamma_i(0)=v_i$. Let $\sphere^{n-1}$ be the equator $\sphere^{n}
\cap(\euclid^{n} \times \{0\})$ and  ${\hat t} = \min\left\{
\dsph(v_i,\sphere^{n-1})\right\}$ be the first time that some
$\gamma_i(t)$ reaches the equator $\sphere^{n-1}$. Denote
$u_i=\gamma_i({\hat t})$. Note that by the construction, the
vertices $u_i\in\sphere^{n-1}$ for $i=1,\ldots,m$ and each $u_{m+i}$
lies in the open hemi-sphere $\sphere^n\cap(\euclid^n\times[-1,0))$
of~$\sphere^n$ for~$i\geq 1$.

As a corollary of Lemma~\ref{lemma-extend}, we have,
\begin{coro}
For $t\in(0,\hat t\,]$, $\dsph(\gamma_i(t),\gamma_j(t)) >
\dsph(v_i,v_j)$ for $i\ne j$.  In particular,
$$
\dsph(u_i,u_j) > \dsph(v_i,v_j) = \pi-\sigma_{ij}.
$$
\end{coro}


\begin{prop} \label{prop-nohalfball}
\begin{enumerate}
\item
The set $\left\{ v_1, \ldots, v_m \right\}$ does not lie in any open
hemi-sphere in $\partial B_s$.
\item
The vectors $u_1, \ldots, u_{n+1}$ do not lie in any open
hemi-sphere in~$\sphere^n$. In particular, the vectors $u_1, ...,
u_{n+1}$ are linearly dependent.
\end{enumerate}
\end{prop}
\begin{proof}
To prove the first statement, we suppose otherwise. Then there is a
unit vector $w$ so that the inner product $(w, v_i) >0$ for
$i=1,..., m$. Now move the center $s$ along the great circle $w_t=
\frac{ (1-t) s + t w}{ ||(1-t) s + tw||}$ where $t \in (0, 1)$. An
easy calculation using $(w, v_i) > 0$ for $i=1,..., m$ shows that
$\dsph(w_t, v_j) < r$ for $t >0$ small and all $j$. This contradicts
the assumption that $B_s$ has the smallest radius.

To see the second statement, suppose otherwise that $u_1, ...,
u_{n+1}$ lie in an open hemi-sphere in $S^n$.  Then the open
hemi-sphere intersects~$\sphere^{n-1}$ in an open hemi-sphere. Since
$u_1, ..., u_m$ are in $\sphere^{n-1}$,  follows that, $u_1, \ldots,
u_{m}$ lie in an open hemi-sphere in~$\sphere^{n-1}$. The spherical
radial rays from~$s$ determine a radial projection between $\partial
B_s$ and~$\sphere^{n-1}$ such that $v_i$'s correspond to~$u_i$'s for
$i=1,\ldots,{n+1}$. Furthermore, the radial projection sends
hemi-spheres to hemi-spheres. Thus, $v_1,\ldots,v_{m}$ also lie in
an open half $(n-1)$-ball in~$\partial B_s$. This contradicts part
(1).

Since any $n+1$ independent unit vectors in $\sphere^n$ lie in an
open hemi-sphere, the last statement follows.

\end{proof}

\begin{proof}[First proof of\/ {\bf M1}]
By proposition 6, there is an $n$-dimensional linear subspace $P$
of~$\euclid^{n+1}$ containing the set $\left\{ u_1, \ldots, u_{n+1}
\right\}$.  Then these points lie in the $(n-1)$-sphere denoted by
$\sphere_1^{n-1} = \sphere^n\cap P$.  By
Proposition~\ref{prop-nohalfball}, $\left\{ u_1, \ldots, u_{n+1}
\right\}$ does not lie in any open hemi-sphere of $\sphere_1^{n-1}$.
Now, we will make use of the following result to finish.
\begin{lemma}
{\bf \cite[Lemma 5]{GuoLuo}}  Let  $\left\{ u_1, \ldots, u_{n+1}
\right\} \subset \sphere^{n-1}$ which does not lie in any open
hemi-sphere of\/ $\sphere^{n-1}$.  For every $\varepsilon >0$, there
is a set $\left\{ w_1, \ldots, w_{n+1} \right\} \subset
\sphere^{n-1}$ such that it does not lie in any hemi-sphere of\/
$\sphere^{n-1}$ and $\dsph(w_i,u_i)<\varepsilon$ for all~$i$.
\end{lemma}
By this lemma, for $\varepsilon = \dfrac{1}{2} \min\left\{
d(u_i,u_j)-d(v_i,v_j) : i\ne j \right\}$, we find the points $w_1,
\ldots, w_{n+1} \in \sphere_1^{n-1}$ such that $d(w_i,w_j)<
\varepsilon$ for all~$i$ and $\left\{ w_1, \ldots, w_{n+1} \right\}$
does not lie in any hemi-sphere in $\sphere_1^{n-1}$. By the choice
of~$\varepsilon$, we have $d(w_i,w_j) > d(v_i,v_i)$ for all~$i\ne
j$.  By Lemma~\ref{lemma-EuclidCpt}, $\sE = \left\{ x\in P : \langle
(x-w_i),w_i\rangle \leq 0\right\}$ is an Euclidean $n$-simplex whose
dihedral angles are given by $\pi-d(w_i,w_j) < \pi-d(v_i,v_i) =
\sigma_{ij}$.

This completes the proof of Theorem~{\bf M1}.
\end{proof}

The geometric relationship between the center $s^*$ of the dual
simplex is very interesting. In fact, due to the convexity, we see
that we always have $s\in\sS^*$. The following two propositions
describe the geometric configuration about the vertices $v_i$'s, the
corresponding $u_i$'s and the center~$s$.

\begin{prop} \label{prop-circumcenter}
The followings are true when $s$~lies in the interior of~$\sS^*$.
\begin{enumerate}
\item
$m = n+1$.
\item
$B_s$ is the $n$-ball circumscribing $\sS^*$, i.e., $v_i \in
\partial B_s$ for all~$i=1,\ldots,n+1$.
\item
The set $\left\{ u_1, \ldots, u_{n+1} \right\}$ does not lie in any
hemi-sphere of $\sphere^{n-1}$.
\end{enumerate}
\end{prop}
\begin{proof}
The first two results follow directly from a special case
($\ell=n+1$) of Lemma~\ref{lemma-smallAtSurf} below.  To get the
last statement, one only needs to follow the argument of
Proposition~\ref{prop-nohalfball}.
\end{proof}
Note that the converse is not true, i.e., even if $m=n+1$, one may
have $s\in\partial\sS^*$.

\begin{prop} \label{prop-bdrycenter}
The followings are true when $s$~lies on the boundary of~$\sS^*$.
\begin{enumerate}
\item
There is an integer~$\ell\leq n$ with $2\leq \ell \leq m \leq n+1$
such that $\ell-1$ is the minimum dimension of a face of~$\sS^*$
which contains~$s$.
\item
$s$~lies in the interior of the face of~$\sS^*$ determined by
$v_1,\ldots,v_\ell$.
\item
$s$~is the center of a geodesic $(\ell-1)$-sphere circumscribing
$\left\{ u_1, \ldots, u_\ell \right\}$.
\item\label{item-bdrycenterface}
$\left\{ u_1, \ldots, u_\ell \right\}$ is the vertex set of a
compact Euclidean $(\ell-1)$-simplex with the origin as
circumcenter.
\item\label{item-bdrycenterDsum}
$\left\{ u_1, \ldots, u_\ell \right\} \subset \euclid^n\times\{0\}$
is of rank~$(\ell-1)$ and $\left\{ u_{\ell+1}, \ldots, u_{n+1}
\right\}$ is linearly independent.  In addition, $\left\{ u_1,
\ldots, u_\ell, \ldots, u_{n+1} \right\}$ is of rank~$n$.
\end{enumerate}
\end{prop}
The following lemma is useful in the proofs of both propositions.
\begin{lemma}\label{lemma-smallAtSurf}
If the center~$s$ of~$B_s$ lies in the interior of the
$(\ell-1)$-face $(v_1,\ldots,v_\ell)$ for some $\ell\leq n+1$, then
$B_s\cap\frak S$ is the $(\ell-1)$-ball circumscribing
$(v_1,\ldots,v_\ell)$, where $\frak S$ is the totally geodesic
$(\ell-1)$-sphere containing $\left\{v_1,\ldots,v_\ell\right\}$.
\end{lemma}
\begin{proof}
It is sufficient to show that $v_1,\ldots,v_\ell\in\partial B_s$. If
$\partial B_s\cap\frak S = \frak S$, then we are done.  If $\partial
B_s\cap\frak S \ne \frak S$, suppose some of $v_i$'s lie in the
interior of~$B_s$ in~$\sphere^n$. Without loss of generality, let
$k<\ell$ and $\left\{v_1,\ldots,v_k\right\}\subset
\partial B_s$ while $v_{k+1}, \ldots, v_\ell \in B_s$.  Since
$s$~lies in the interior of $(v_1,\ldots,v_\ell)$ and
$\text{radius}(B_s)<\pi/2$, it does not lie in the geodesic
$(k-1)$-sphere spanned by $v_1,\ldots,v_k$. By the proof of
Proposition~\ref{prop-nohalfball}, we may perturb~$s$ to~$s'$ and
have a ball of smaller radius.
\end{proof}

\begin{proof}[Proof of Proposition~\ref{prop-bdrycenter}]
Let $\ell-1$ be the lowest dimension of a face of
$(v_1,\ldots,v_{n+1})$ that contains the center~$s$.  Obviously,
$\ell\geq 2$ and by the minimality of~$\ell$, $s$~lies in the
interior of the face. Without loss of generality, assume this face
has vertices $\left\{v_1,\ldots,v_\ell\right\}$ and it determines a
totally geodesic $(\ell-1)$-sphere~$\frak S$.  By
Lemma~\ref{lemma-smallAtSurf}, $B_s\cap\frak S$ is the
$(\ell-1)$-ball circumscribing $\left\{v_1,\ldots,v_\ell\right\}$.
Thus, $\ell \leq m$.  Using the same argument as in
Proposition~\ref{prop-nohalfball}, we can see that
$\left\{u_1,\ldots,u_\ell\right\}$ does not lie in any open half
$(\ell-1)$-ball of $\frak S\cap\sphere^{n-1}$.  Thus, it determines
a compact Euclidean $(\ell-1)$-simplex.  The last statement now
follows from the nondegeneracy of~$\sS^*$ and a dimension count.
\end{proof}

Based on Propositions~\ref{prop-circumcenter}
and~\ref{prop-bdrycenter}, we are giving a more explicit alternative
proof for Theorem~{\bf M1}.
\begin{proof}[Second proof of\/ {\bf M1}]
First, let us consider the case that $s\in\left(\sS^*\right)^\circ$.
By Proposition~\ref{prop-circumcenter},  $B_s$~is the circumscribe
$n$-ball of~$\sS^*$ and for all $i,j$, we have
$$
\dsph(u_i,u_j) > \dsph(v_i,v_j) = \pi - \sigma_{ij}.
$$
Moreover,
$$
\left\{ u_1, \ldots, u_{n+1}\right\} \subset \sphere^{n-1} \subset
\euclid^n\times\{0\} \subset \euclid^{n+1};
$$
but it does not lie in any closed half $(n-1)$-ball of
$\sphere^{n-1}$.

Let $\sE$ be the subset of~$\euclid^n\times\{0\}$ bounded by the
codimension-1 hyperplanes tangent to $\sphere^{n-1}$ at the $u_i$'s.
Since the $u_i$'s do not lie in any closed half-space, these tangent
hyperplanes bound a compact Euclidean $n$-simplex $\sE$
in~$\euclid^n\times\{0\}$ with dihedral angles $\xi_{ij} = \pi -
\dsph(u_i,u_j) < \sigma_{ij}$.  So, $\sE$ is the required Euclidean
$n$-simplex.

In the case that $s\in\partial\sS^*$, by
Proposition~\ref{prop-bdrycenter},
statement~(\ref{item-bdrycenterface}), there exists $a_i>0$,
$i=1,\ldots,\ell$, such that $\displaystyle \sum_{i=1}^{\ell}
a_i\,u_i = 0.$ Take arbitrarily small $\delta>0$ and let
$$
w_i = \begin{cases}
\quad u_i - \delta\left(u_{\ell+1}+\cdots+u_{n+1}\right) &\quad i=1,\ldots, \ell, \\
\quad u_i &\quad i=\ell+1,\ldots,n+1\,.\end{cases}
$$
One may choose $b_i>0$ as follows,
$$
b_i = \begin{cases} \quad
a_i\left/\left(\sum_{q=1}^{\ell}a_q\right)\right. &\quad
i=1,\ldots,\ell \,;\\
\quad\delta &\quad i=\ell+1,\ldots, n+1\,.
\end{cases}
$$
Then,
$$
\sum_{i=1}^{n+1} b_i w_i = \sum_{i=1}^{\ell}
\frac{a_i}{\sum_{q=1}^{\ell}a_q} u_i \,-\, \sum_{i=1}^{\ell}
\frac{a_i\,\delta}{\sum_{q=1}^{\ell}a_q}\sum_{j=\ell+1}^{n+1} u_j
\,+\,\delta \sum_{i=\ell+1}^{n+1} u_i = 0\,.
$$
Next, we will prove that one may choose $\delta>0$ such that any
subset of $n$~vectors among $\left\{w_1, \ldots,\ldots,
w_{n+1}\right\}$ is linearly independent.  We will consider the
subset $\left\{w_1, \ldots,\crossout{w}_q,\ldots, w_{n+1}\right\}$
in the cases that $q\leq \ell$ or $q\geq \ell+1$.

Let $q\leq \ell$ and $\displaystyle\sum_{q\ne i=1}^{n+1} x_iw_i =
0$. Substituting the expressions of $w_i$'s, we have
$$
\sum_{\substack{i=1\\ i\ne q}}^{\ell} x_iu_i + \sum_{i=\ell+1}^{n+1}
\left( x_i - \delta\sum_{\substack{j=1\\ j\ne q}}^{\ell}x_j
\right)u_i = 0\,.
$$
Observe that if $q\leq \ell$, by (\ref{item-bdrycenterDsum}) of
Proposition~\ref{prop-bdrycenter},
$\left\{u_{1},\ldots,\crossout{u}_q,\ldots,u_{n+1} \right\}$ is
linearly independent.  The above equation implies that $x_i=0$ for
all $i\ne q$.

In the case that $q\geq \ell+1$ and $\displaystyle\sum_{q\ne
i=1}^{n+1} x_iw_i = 0$ for some $x_i$'s and a certain $\delta>0$. We
claim that only one specific $\delta$ may have nontrivial $x_i$'s.
By substituting the expressions of $w_i$'s, we have
\begin{equation}\tag{$\star$}
\sum_{i=1}^{\ell}x_iu_i - \delta \left(\sum_{j=1}^{\ell}
x_j\right)u_q + \sum_{\substack{i=\ell+1\\ i\ne q}}^{n+1} \left( x_i
- \delta\sum_{j=1}^{\ell} x_j \right)u_i = 0\,.
\end{equation}
Since $\left\{u_1,\ldots,u_{n+1}\right\}$ has rank~$n$, the above
equation has a one-dimensional space for the coefficients. If there
are $\delta_1, \delta_2 > 0$ and corresponding $x_i^{(1)}$,
$x_i^{(2)}$ which satisfy the above equation~($\star$), one can
conclude that
$$
\delta_1 = \delta_2 \quad \text{or} \quad \sum_{i=1}^{\ell}
x_i^{(1)} = \sum_{i=1}^{\ell} x_i^{(2)} = 0.
$$
We will rule out the second alternative.  Suppose there is a
non-trivial set of $x_i$'s with $\sum_{i=1}^\ell x_i=0$ such that
$(\star)$ holds.  Then, equation~($\star$) becomes
$$
\sum_{i=1}^{\ell} x_iu_i + \sum_{\substack{i=\ell+1\\ i\ne q}}^{n+1}
x_i u_i = 0.
$$
By (\ref{item-bdrycenterDsum}) of Proposition~\ref{prop-bdrycenter},
the vectors $\left\{ u_i \right\}_{i=1}^\ell$ and $\left\{ u_i
\right\}_{i=\ell+1}^{n+1}$ span direct summands.  Thus, we must have
simultaneously
$$
\sum_{i=1}^{\ell} x_iu_i = 0, \qquad \sum_{\substack{i=\ell+1\\ i\ne
q}}^{n+1} x_i u_i =0.
$$
However, $\displaystyle\sum_{i=1}^{\ell} x_iu_i = 0$ together with
$\displaystyle\sum_{i=1}^{\ell} x_i = 0$ contradict that $u_{1},
\ldots, u_{\ell}$ form a compact Euclidean simplex.  Consequently,
one must have $\delta_1 = \delta_2$.

Thus, by \cite[Lemma~4]{GuoLuo}, there is sufficiently small
$\delta>0$ such that the vertices $w_i$,~$i=1,\ldots,n+1$ span an
$n$-dimensional space~$L\subset\R^{n+1}$ and they define a compact
Euclidean $n$-simplex in ~$L$. Furthermore, $\norm{w_i-u_i}$ can be
made arbitrarily small. Let $\sE$ be the Euclidean $n$-simplex
in~$L$ dual to $w_i$'s.  In other words, if $w_i$'s are normalized,
$\sE$ is bounded by the tangent hyperplanes to $\sphere^n\cap L$
at~$w_i$. Its dihedral angles $\xi_{ij}$ satisfy that
$\left|\xi_{ij}-(\pi-\ell_{ij})\right|<\varepsilon$ for arbitrarily
small~$\varepsilon>0$. Hence,
$$
\xi_{ij} < \pi-\ell_{ij}+\varepsilon < \pi-\dsph(v_i,v_j) =
\sigma_{ij}.
$$
This completes the proof of the theorem.
\end{proof}


\end{document}